\documentclass[a4paper,12pt]{article}
\usepackage{amsmath, amssymb, amsthm}
\usepackage{graphicx} 

\begin{document}

\def\Homeo{{\rm Homeo}}\def\MM{{\mathcal M}}\def\Z{{\mathbb Z}}
\def\BB{{\mathcal B}}\def\Sym{{\rm Sym}}\def\Stab{{\rm Stab}}
\def\Ad{{\rm Ad}}\def\Im{{\rm Im}}\def\N{{\mathbb N}}
\def\cycl{{\rm cycl}}\def\Supp{{\rm Supp}}\def\Id{{\rm Id}}
\def\sg{{\rm sg}}\def\S{{\mathbb S}}

%%%%%%%%%%%%%%%%%%%%%%%%%%%%%%%%%%%%%%%%%%%%%%%%%%%%%%%%%%%%%
\title{\bf{Small index subgroups of the mapping class group}}
 
\author{\textsc{Luis Paris}}

\date{\today}

\maketitle

\begin{abstract}
\noindent
We prove that the mapping class group of a closed oriented surface of genus $\rho \ge 3$ has no 
proper subgroup of index $\le 4 \rho +4$.
\end{abstract}


%%%%%%%%%%%%%%%%%%%%%%%%%%%%%%%%%%%%%%%%%%%%%%%%%%%%%%%%%%%%%%%%%%%%%%%%%%%%%%%%%%%%%%%%%
\vskip3truecm\noindent
Let $\Sigma$ be a closed oriented surface of genus $\rho \ge 1$, and let $\Homeo^+ (\Sigma)$ be the 
group of homeomorphisms $h: \Sigma \to \Sigma$ that preserve the orientation. The {\it mapping 
class group} of $\Sigma$ is defined to be the group $\MM (\Sigma) = \pi_0 (\Homeo^+ (\Sigma))$ of isotopy classes 
of elements of $\Homeo^+(\Sigma)$.
This group plays an important role in low-dimensional topology, in Teichm\"uller theory, and in group
theory. We refer to \cite{Ivano1} for a general survey on the subject, and turn to the main object
of the present paper: the finite index subgroups of $\MM (\Sigma)$.

\bigskip\noindent
A group $G$ is called {\it residually finite} if, for all $g \in G \setminus \{ 1 \}$, there exists a 
homomorphism $\varphi: G \to H$ such that $H$ is finite and $\varphi (g) \neq 1$. By \cite{Gross1}, the 
mapping class group $\MM (\Sigma)$ is residually finite. This means that $\MM (\Sigma)$ has many finite 
index subgroups, numerous enough to separate its elements. Nevertheless, we are going to prove:

\bigskip\noindent
{\bf Theorem 1.} {\it Let $\Sigma$ be a closed oriented surface of genus $\rho \ge 3$. Then $\MM 
(\Sigma)$ has no proper subgroup of index $\le 4 \rho +4$.}

\bigskip\noindent
{\bf Remark.} If $\rho =1$, then $H_1 (\MM (\Sigma)) \simeq \Z/ 12\Z$, thus $\MM (\Sigma)$ has (normal) 
subgroups of index 2 and 3. If $\rho=2$, then $H_1 (\MM (\Sigma)) \simeq \Z/10\Z$ (see \cite{Mumfo1}), 
thus $\MM (\Sigma)$ has (normal) subgroups of index 2 and 5. So, the hypothesis $\rho \ge 3$ is needed 
in Theorem 1. Note also that the statement is valid for all the subgroups, not only for the normal ones.

\bigskip\noindent
The rest of the article is dedicated to the proof of Theorem 1. So, from now on, $\Sigma$ denotes a 
closed oriented surface of genus $\rho \ge 3$.

\bigskip\noindent
The two main tools that we use in our proof are the braid groups and the Dehn twists. For $n \ge 2$, we 
denote by $\BB_n$ the braid group on $n$ strands. This has a presentation with generators $\sigma_1, 
\dots, \sigma_{n-1}$ and relations
\[
\begin{array}{cl}
\sigma_i \sigma_j = \sigma_j \sigma_i&\quad \text{if } |i-j| \ge 2\,,\\
\sigma_i \sigma_j \sigma_i = \sigma_j \sigma_i \sigma_j &\quad \text{if } |i-j|=1\,.
\end{array}
\]
The Dehn twist along an essential circle $a: \S^1 \hookrightarrow \Sigma$ will be denoted by $\tau_a$. 
We assume the reader to be familiar with these two concepts, and we refer to \cite{Birma1} for a 
detailed text on them.

\bigskip\noindent
The following is easy to prove and can be found, for instance, in \cite{Birma2}.

\bigskip\noindent
{\bf Lemma 2.} {\it Let $a,b$ be two essential circles which intersect transversely. Then
\[
\begin{array}{cl}
\tau_a \tau_b = \tau_b \tau_a &\quad\text{if } a \cap b = \emptyset\,,\\
\tau_a \tau_b \tau_a = \tau_b \tau_a \tau_b &\quad \text{if } |a \cap b| =1\,.
\end{array}
\]}

\bigskip\noindent
Now, consider the essential circles $a_0, a_1, \dots, a_{2 \rho+1}$ drawn in Figure 1, and, for $0 \le 
i\le 2 \rho +1$, let $\tau_i$ denote the Dehn twist along $a_i$.

%%%%%%%%%%%
\begin{figure}[htb]
\bigskip
\centerline{
\setlength{\unitlength}{.5cm}
\begin{picture}(20,4)
\put(0,0){\includegraphics[width=10cm]{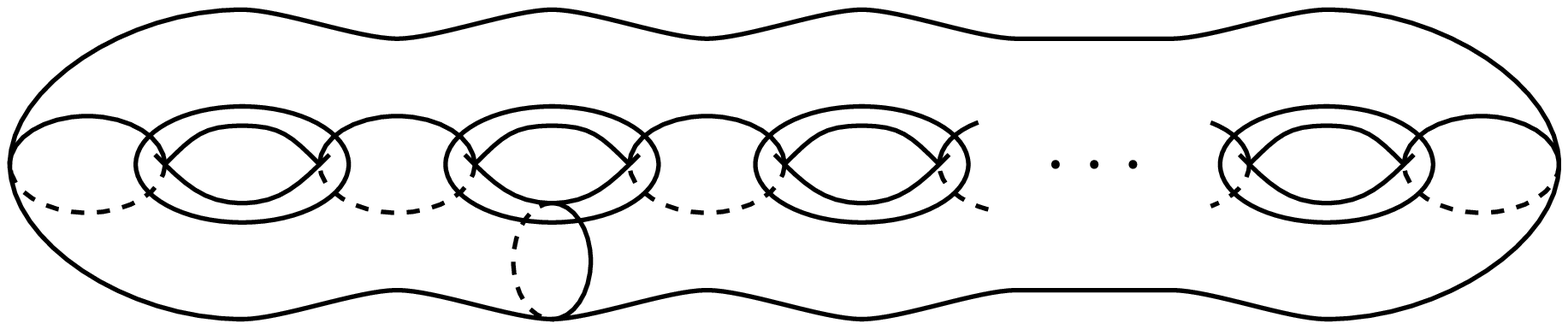}}
\put(0.8,2.9){\small $a_1$}
\put(2.5,3){\small $a_2$}
\put(4.7,2.9){\small $a_3$}
\put(6.5,3){\small $a_4$}
\put(8.7,2.9){\small $a_5$}
\put(10.5,3){\small $a_6$}
\put(16.5,3){\small $a_{2\rho}$}
\put(18.2,2.9){\small $a_{2 \rho +1}$}
\put(7.6,0.6){\small $a_0$}
\end{picture}}
\bigskip
\centerline{{\bf Figure 1.} Essential circles in $\Sigma$.}
\end{figure}
%%%%%%%%%

\bigskip\noindent
{\bf Proposition 3.} {\it
\begin{enumerate}
\item
The mapping $\sigma_i \mapsto \tau_i$, $1 \le i\le 2 \rho +1$, induces a homomorphism $\iota: \BB_{2 
\rho +2} \to \MM (\Sigma)$. Furthermore, we have the relations
\[
\begin{array}{cl}
\tau_0 \tau_i = \tau_i \tau_0 &\quad \text{for } 1 \le i\le 2 \rho +1 \text{ and } i \neq 4\,,\\
\tau_0 \tau_4 \tau_0 = \tau_4 \tau_0 \tau_4\,.
\end{array}
\]
\item
{\rm (Humphries \cite{Humph1})} The Dehn twists $\tau_0, \tau_1, \dots, \tau_{2 \rho}$ generate $\MM 
(\Sigma)$.
\item
{\rm (Powell \cite{Powel1})} $H_1 (\MM (\Sigma))= \{ 0\}$.
\end{enumerate}}

\bigskip\noindent
{\bf Remark.} Proposition 3.1 is a straightforward consequence of Lemma 2. Note also that $\iota: 
\BB_{2 \rho+2} \to \MM (\Sigma)$ is not injective (see \cite{LabPar1}).

\bigskip\noindent
For $k \ge 1$, we denote by $\Sym_k$ the symmetric group of $\{ 1,2, \dots, k\}$. Let $G$ be a group. A 
homomorphism $\varphi: G \to \Sym_k$ is called {\it transitive} if $\Im\, \varphi$ acts transitively on 
$\{ 1,2, \dots, k\}$. If $\varphi: G \to \Sym_k$ is a transitive homomorphism, then $\Stab_\varphi (1) 
= \{ g \in G; \varphi(g)(1) =1 \}$ is an index $k$ subgroup of $G$. Reciprocally, if $H$ is an index 
$k$ subgroup of $G$, then there exists a transitive homomorphism $\varphi: G \to \Sym_k$ such that $H= 
\Stab_\varphi(1)$ (take the action of $G$ on the left cosets of $H$). So, in order to prove Theorem 1, 
it suffices to prove that there is no transitive homomorphism $\varphi: \MM (\Sigma) \to \Sym_k$ for 
any $2 \le k\le 2 \rho +4$.

\bigskip\noindent
For $u \in Sym_k$, we denote by $\Ad_u: \Sym_k \to \Sym_k$ the inner automorphism which sends $v$ to $u 
v u^{-1}$ for all $v \in \Sym_k$. We say that two homomorphisms $\varphi, \varphi' : G \to \Sym_k$ are 
{\it conjugate} if there exists $u \in \Sym_k$ such that $\varphi'= \Ad_u \circ \varphi$. To be 
conjugate is an equivalence relation that we denote by $\varphi \sim \varphi'$. The {\it product} of 
two homomorphisms $\varphi_1: G \to \Sym_{k_1}$ and $\varphi_2: G \to \Sym_{k_2}$ is the homomorphism 
$\varphi_1 \times \varphi_2 : G \to \Sym_{k_1+k_2}$ defined by
\[
(\varphi_1 \times \varphi_2) (g)(i) = \left\{
\begin{array}{cl}
\varphi_1 (g)(i) &\quad \text{if } 1 \le i\le k\,,\\
\varphi_2(g)(i-k_1)+k_1 &\quad\text{if } k_1+1 \le i\le k_1+k_2\,.
\end{array}\right.
\]
Clearly, for every homomorphism $\varphi: G \to \Sym_k$, there exists a partition $(k_1, k_2, \dots, 
k_p)$ of $k$ and a collection $\{ \varphi_j: G \to \Sym_{k_j}; 1 \le j\le p\}$ of transitive 
homomorphisms such that
\[
\varphi \sim \varphi_1 \times \varphi_2 \times \cdots \times \varphi_p\,.
\]
Moreover, this decomposition is unique up to a permutation of the factors.

\bigskip\noindent
Using the above defined product of homomorphisms, it is easily shown that the existence of a transitive 
homomorphism $\varphi: \MM (\Sigma) \to \Sym_k$ for some $2 \le k\le 4 \rho +4$ is equivalent to the 
existence of a non-trivial homomorphism $\varphi: \MM (\Sigma) \to \Sym_{4 \rho +4}$.

\bigskip\noindent
We set $n= 2 \rho +2$ and we turn now to give the classification of the homomorphisms from $\BB_n$ to 
$\Sym_{2n}$ (see Corollary 7). Then, given a homomorphism $\varphi: \MM (\Sigma) \to \Sym_{2n} = 
\Sym_{4 \rho +4}$, our strategy for proving that $\varphi$ is trivial will consist on exploring the 
different possibilities for the composition $\psi= \varphi \circ \iota: \BB_n \to \Sym_{2n}$, where 
$\iota: \BB_n = \BB_{2 \rho +2} \to \MM (\Sigma)$ is the homomorphism of Proposition 3.1.

\bigskip\noindent
Let $k \in \N$ 
and let $w \in \Sym_k$. Then the mapping $\sigma_i \mapsto w$, $1 \le i\le n-1$, induces a homomorphism 
$\cycl_w: \BB_n \to \Sym_k$ called a {\it cyclic homomorphism}. Note that $\cycl_w$ is transitive if 
and only if $w$ is a cycle of length $k$. On the other hand, the mapping $\sigma_i \mapsto (i,i+1) \in 
\Sym_n$, $1 \le i\le n-1$, induces a transitive homomorphism $\rho_S: \BB_n \to \Sym_n$ called the {\it standard 
homomorphism}.

\bigskip\noindent
{\bf Theorem 4} (Artin \cite{Artin1}). {\it Let $n \ge 7$.
\begin{enumerate}
\item
For $k <n$, every transitive homomorphism $\psi: \BB_n \to \Sym_k$ is cyclic.
\item
Every transitive homomorphism $\psi: \BB_n \to \Sym_n$ is either cyclic or conjugate to the standard 
homomorphism $\rho_S$.
\end{enumerate}}

\bigskip\noindent
{\bf Remark.} For $n=4$ or $6$, there are transitive homomorphisms $\BB_n \to \Sym_n$ which are neither 
cyclic, nor conjugate to $\rho_S$: there are 3 for $n=4$, and 1 for $n=6$, up to conjugation (see 
\cite{Artin1}). However, in this paper we are concerned only by the case $n \ge 8$.

\bigskip\noindent
{\bf Proposition 5} (Lin \cite{Lin1}). {\it Let $n \ge 7$.
\begin{enumerate}
\item
The mapping
\[
\sigma_i \mapsto (i,n+i+1, n+i, i+1)\,, \quad 1 \le i\le n-1\,,
\]
induces a transitive homomorphism $\rho_{L\,1}: \BB_n \to \Sym_{2n}$.
\item
The mapping
\[
\sigma_i \mapsto \prod_{j=1}^{i-1} (j,n+j) \cdot (i,n+i+1,n+i,i+1) \cdot \prod_{j=i+2}^n (j,n+j)\,, 
\quad 1 \le i\le n-1\,,
\]
induces a transitive homomorphism $\rho_{L\,2}: \BB_n \to \Sym_{2n}$.
\item
The mapping 
\[
\sigma_i \mapsto \prod_{j=1}^{i-1} (j,n+j) \cdot (i,n+i+1) (i+1,n+i) \cdot \prod_{j=i+2}^n (j,n+j)\,, 
\quad 1 \le i\le n-1\,,
\]
induces a transitive homomorphism $\rho_{L\,3}: \BB_n \to \Sym_{2n}$.
\end{enumerate}}

\bigskip\noindent
The homomorphisms $\rho_{L\,1}, \rho_{L\,2}, \rho_{L\,3}$ of Proposition 5 will be called {\it Lin 
homomorphisms}.

\bigskip\noindent
{\bf Theorem 6} (Lin \cite{Lin1}). {\it Let $n \ge 7$.
\begin{enumerate}
\item
For $n+1 \le k \le 2n-1$, every transitive homomorphism $\psi: \BB_n \to \Sym_k$ is cyclic.
\item
Every transitive homomorphism $\psi: \BB_n \to \Sym_{2n}$ is either cyclic or conjugate to some 
$\rho_{L\,j}$, where $j \in \{1,2,3\}$.
\end{enumerate}}

\bigskip\noindent
{\bf Corollary 7.} {\it Let $n \ge 7$. Let $\psi: \BB_n \to \Sym_{2n}$ be a homomorphism. Then either
\begin{itemize}
\item[(a)]
$\psi$ is cyclic; or
\item[(b)]
$\psi \sim \rho_S \times \cycl_w$ for some $w \in \Sym_n$; or
\item[(c)]
$\psi \sim \rho_S \times \rho_S$; or
\item[(d)]
$\psi \sim \rho_{L\,j}$ for some $j \in \{1,2,3\}$.
\end{itemize}}

\bigskip\noindent
{\bf Proof of Theorem 1.} Set $n = 2 \rho +2 \ge 8$. Let $\varphi: \MM (\Sigma) \to \Sym_{4 \rho +4} = 
\Sym_{2n}$ be a homomorphism. As pointed out before, we must prove that $\varphi$ is trivial. We set 
$\psi = \varphi \circ \iota: \BB_n \to \Sym_{2n}$, where $\iota: \BB_n \to \MM (\Sigma)$ is the 
homomorphism of Proposition 3.1. By Corollary 7, we have either
\begin{itemize}
\item[(a)]
$\psi=\cycl_w$ for some $w \in \Sym_{2n}$; or
\item[(b)]
$\psi \sim \rho_S \times \cycl_w$ for some $w \in \Sym_n$; or
\item[(c)]
$\psi \sim \rho_S \times \rho_S$; or
\item[(d)]
$\psi \sim \rho_{L\,j}$ for some $j \in \{1,2,3\}$.
\end{itemize}

\bigskip\noindent
{\bf Case 1:} Suppose $\psi= \cycl_w$, where $w \in \Sym_{2n}$.

\bigskip\noindent
Set $f = \tau_4 \tau_0$ and observe that $f \tau_4 f^{-1} = \tau_0$, and $f \tau_1 f^{-1} = \tau_1$. 
We have $\varphi (\tau_i) = \psi (\sigma_i) = w$ for all $1 \le i \le 2\rho$. Moreover,
\[
\varphi (\tau_0) = \varphi(f) \varphi(\tau_4) \varphi(f^{-1}) = \varphi(f) \varphi(\tau_1) \varphi(f^{-
1}) = \varphi(\tau_1) = w\,.
\]
Since $\MM (\Sigma)$ is generated by $\tau_0, \tau_1, \dots, \tau_{2\rho}$ (see Proposition 3.2), it 
follows that the image of $\varphi$ is the cyclic group generated by $w$: an abelian group. Since $H_1 
(\MM (\Sigma)) = \{ 0\}$ (see Proposition 3.3), we conclude that $\Im\, \varphi$ is trivial, that is, 
$w=1$ and $\varphi$ is trivial.

\bigskip\noindent
{\bf Case 2:} Suppose $\psi= \rho_S \times \cycl_w$, where $w \in \Sym_n$.

\bigskip\noindent
For $v \in \Sym_{2n}$, we define the {\it support} of $v$ to be $\Supp(v) = \{ i \in \{ 1, \dots, 2n\}; 
v(i) \neq i \}$. We will often use the fact that two elements commute if their supports are disjoint.

\bigskip\noindent
Let $w_0 \in \Sym_{2n}$ be defined by
\[
w_0(i)= \left\{
\begin{array}{cl}
i &\quad\text{if } 1 \le i\le n\,,\\
w(i-n)+n &\quad\text{if } n+1\le i\le 2n\,.
\end{array}\right.
\]
So, $\Supp (w_0) \subset \{ n+1, n+2, \dots, 2n\}$, and
\[
\varphi (\tau_i) = \psi( \sigma_i) = (i,i+1) w_0\,, \quad \text{for all } 1 \le i\le n-1\,.
\]
Set $u_0 = \varphi(\tau_0)$. The permutation $u_0$ commutes with $\varphi (\tau_1) \varphi (\tau_2)^{-
1} = (1,2,3)$, with $\varphi (\tau_2) \varphi(\tau_3)^{-1} = (2,3,4)$, and with $\varphi (\tau_1) 
\varphi (\tau_3)^{-1} = (1,2) (3,4)$, thus
\[
u_0 (\{ 1,2,3\}) = \{1,2,3\}\,, \quad u_0 (\{2,3,4\}) = \{ 2,3,4\}\,,
\]
and either
\[ \begin{array}{cl}
u_0 (\{1,2\}) = \{1,2\} \quad \text{and} \quad u_0 (\{3,4\}) = \{3,4\}\,, &\quad \text{or}\\
u_0 (\{1,2\}) = \{3,4\} \quad \text{and} \quad u_0 (\{3,4\}) = \{1,2\}\,.
\end{array}\]
This clearly implies that $u_0(i)=i$ for all $1 \le i\le 4$. Similarly, the fact that $u_0$ commutes 
with $\varphi (\tau_i) = (i,i+1) w_0$ for all $5 \le i\le n-1$ implies that $u_0(i)=i$ for all $5 \le 
i\le n$. So, $\Supp (u_0) \subset \{ n+1, \dots, 2n \}$. Now, from the equality 
$\tau_0 \tau_4 \tau_0 = \tau_4 \tau_0 \tau_4$ follows
\[
(4,5) u_0 w_0 u_0 = (4,5)(4,5) w_0 u_0 w_0 = w_0 u_0 w_0\,,
\]
and this last equality cannot hold: a contradiction.

\bigskip\noindent
{\bf Case 3:} Suppose $\psi= \rho_S \times \rho_S$.

\bigskip\noindent
We have
\[
\varphi (\tau_i) = \psi (\sigma_i) = (i,i+1) (n+i,n+i+1)\,, \quad \text{for all } 1 \le i\le n-1\,.
\]
Set $u_0= \varphi(\tau_0)$. Recall that $\tau_0 = f \tau_4 f^{-1}$, where $f = \tau_4 \tau_0$, thus 
$u_0$ is conjugate to $\varphi (\tau_4) = (4,5) (n+4,n+5)$, therefore $u_0$ is the disjoint product of 
two transpositions. Since $u_0$ commutes with $\varphi (\tau_i) = (i,i+1) (n+i,n+i+1)$ for all $1 \le 
i\le 3$, we have
\[
u_0 (\{i,i+1,n+i,n+i+1\}) = \{i,i+1,n+i,n+i+1\}
\]
for all $1 \le i\le 3$, thus $u_0( \{i,n+i\}) = \{i,n+i\}$ for all $1 \le i \le 4$. If $u_0(i) = n+i$ 
then $u_0(i+1) = n+i+1$, since $u_0$ commutes with $\varphi (\tau_i) = (i,i+1) (n+i,n+i+1)$. Similarly, 
if $u_0(i)=i$, then $u_0(i+1) = i+1$. So, either $u_0(i)=i$ and $u_0(n+i)=n+i$ for all $1 \le i\le 4$, 
or $u_0(i)=n+i$ and $u_0(n+i) = i$ for all $1 \le i\le 4$. The latter case cannot hold because $u_0$ is 
the disjoint product of two transpositions, thus $u_0(i)=i$ and $u_0(n+i)=n+i$ for all $1 \le i\le 4$. 
But, the fact that $u_0$ commutes with $\varphi (\tau_i) = (i,i+1) (n+i,n+i+1)$ for all $5 \le i\le n-
1$ also implies that $u_0(i)=i$ and $u_0(n+i)=n+i$ for all $5 \le i\le n$. So, $u_0=\Id$: a 
contradiction.

\bigskip\noindent
{\bf Case 4:} Suppose $\psi= \rho_{L\,1}$.

\bigskip\noindent
Let $\sg: \Sym_{2n} \to \{ \pm 1\}$ denote the signature. Observe that
\[
(\sg \circ \varphi) (\tau_1) = (\sg \circ \psi) (\sigma_1) = \sg((1,n+2,n+1,2))= -1\,,
\]
thus $\sg \circ \varphi : \MM (\Sigma) \to \{ \pm 1\}$ is surjective. This contradicts the fact that 
$H_1 (\MM (\Sigma)) = \{0\}$ (see Proposition 3.3).

\bigskip\noindent
{\bf Case 5:} Suppose $\psi= \rho_{L\,2}$.

\bigskip\noindent
Note that $n=2 \rho +2$ is even, thus
\[
(\sg \circ \varphi) (\tau_1) = (\sg \circ \psi) (\sigma_1) = \sg \left( (1,n+2,n+1,2) \cdot 
\prod_{j=3}^n (j,n+j) \right) = (-1)^{n-1} = -1\,,
\]
therefore $\sg \circ \varphi: \MM (\Sigma) \to \{ \pm 1\}$ is surjective. As pointed out in the 
previous case, such a homomorphism cannot exist because $H_1 (\MM (\Sigma)) = \{ 0\}$.

\bigskip\noindent
{\bf Case 6:} Suppose $\psi= \rho_{L\,3}$.

\bigskip\noindent
We have
\[
\varphi (\tau_i) = \psi (\sigma_i) = \prod_{j=1}^{i-1} (j,n+j) \cdot (i,n+i+1) (i+1, n+i) \cdot 
\prod_{j=i+2}^n (j,n+j)\,,
\]
for all $1 \le i\le n-1$. Set $u_0 = \varphi (\tau_0)$. Recall that $\tau_0$ is conjugate to $\tau_4$, 
thus $u_0$ is conjugate to $\varphi(\tau_4)$, therefore $u_0$ is the disjoint product of $n$ 
transpositions. The permutation $u_0$ commutes with the following three:
\[
\begin{array}{rcl}
\varphi (\tau_1) \varphi (\tau_2) &=& (1,2,3) (n+1,n+2,n+3)\,,\\
\varphi (\tau_1) \varphi (\tau_3) &=& (1,2) (3,4) (n+1,n+2) (n+3,n+4)\,,\\
\varphi(\tau_2) \varphi (\tau_3) &=& (2,3,4) (n+2,n+3,n+4)\,,
\end{array}
\]
thus
\[
\begin{array}{rcl}
u_0( \{ 1,2,3, n+1,n+2,n+3\}) &=& \{ 1,2,3, n+1,n+2,n+3\}\,,\\
u_0(\{ 1,2,3,4, n+1,n+2,n+3,n+4 \}) &=& \{ 1,2,3,4, n+1,n+2,n+3,n+4 \}\,,\\
u_0( \{ 2,3,4, n+2,n+3,n+4 \}) &=& \{ 2,3,4, n+2,n+3,n+4 \}\,,
\end{array}
\]
hence
\[
u_0 (\{ 1,n+1 \}) = \{ 1, n+1 \} \quad \text{and} \quad u_0 (\{ 4,n+4 \}) = \{4,n+4\}\,.
\]
Moreover, the above two equalities plus the fact that $u_0$ commutes with $\varphi (\tau_1) 
\varphi (\tau_3) = (1,2) (3,4) (n+1,n+2) (n+3,n+4)$ imply that we also have
\[
u_0 (\{ 2,n+2 \}) = \{ 2, n+2 \} \quad \text{and} \quad u_0 (\{ 3,n+3 \}) = \{3,n+3\}\,.
\]
Similarly, the fact that $u_0$ commutes with $\varphi (\tau_i)$ for all $5 \le i\le n-1$ implies that 
$u_0 (\{ i,n+i \}) = \{i,n+i\}$ for all $5 \le i\le n$. Finally, since $u_0$ is the disjoint product of 
$n$ transpositions, we conclude that
\[
u_0 = \prod_{j=1}^n (j,n+j)\,.
\]
Now, a direct calculation shows that
\[
u_0 \varphi (\tau_4) u_0 = \varphi (\tau_4) \neq u_0 = \varphi (\tau_4) u_0 \varphi (\tau_4)\,,
\]
and this contradicts the equality $\tau_0 \tau_4 \tau_0 = \tau_4 \tau_0 \tau_4$ of Proposition 3.1.
\qed

%%%%%%%%%%%%%%%%%%%%%%%%%%

%%%%%%%%%%%%%%%%%%%%%%%%%%%%%%%%%%%%%%%%%%%%%%%%%%%%%%%%%%%%%%%%%%%%%%%%%%%%%%%%%%%%%%%%%
\bigskip\bigskip\noindent
{\bf Luis Paris},

\smallskip\noindent 
Institut de Math\'ematiques de Bourgogne, UMR 5584 du CNRS, Universit\'e de Bourgogne, B.P. 
47870, 21078 Dijon cedex, France

\smallskip\noindent
E-mail: {\tt lparis@u-bourgogne.fr}

%%%%%%%%%%%%%%
\end{document}